\documentclass[11pt]{{article}}
\usepackage{CJK}
\usepackage{amsmath}
\usepackage{amsfonts}
\usepackage{mathrsfs}
\usepackage{amssymb,amsmath}
\usepackage{amstext}
\usepackage{array}
\usepackage{theorem}
\usepackage{latexsym}
\usepackage[all]{xy}
\usepackage{amscd}
\usepackage{bbding}
\usepackage{graphics}
\usepackage{graphicx}
\usepackage{bbm}
\usepackage{color}
\usepackage{pgf,pgfarrows,pgfnodes,pgfautomata,pgfheaps}
  
\usepackage[all]{xy}
\def\mod{\mathop{\rm mod}\nolimits}

\def\id{\mathop{\rm id}\nolimits}
\def\pd{\mathop{\rm pd}\nolimits}
\def\gldim{\mathop{\rm gl.dim}\nolimits}

\title{\Large \bf $\tau$-rigid modules for algebras
with radical square zero
}
\author{ Xiaojin Zhang\thanks{{\it E-mail address}:
xjzhang@nuist.edu.cn,\ \
xiaojinzhang.cn@gmail.com}\\
{\footnotesize \it School of Mathematics and Statistics, Nanjing
University of Information Science and Technology},\\{\footnotesize
\it Nanjing 210044, P. R. China}}
\date{}

\begin{document}

\baselineskip=18pt
 \maketitle
 \begin{abstract} In this paper, we show that for an algebra $\Lambda$
 with radical square zero and an indecomposable $\Lambda$-module $M$ such that $\Lambda$ is Gorenstein of finite type or $\tau M$ is $\tau$-rigid,
 $M$ is $\tau$-rigid if and only if the first two projective terms of a minimal projective resolution of $M$
have no on-zero direct summands in common.
 We also determined all $\tau$-tilting modules for Nakayama algebras with
 radical square zero. Moreover, by giving a construction theorem we
 show that a basic connected radical square zero algebra admitting a unique
 $\tau$-tilting module is local.
\end{abstract}

\vspace{0.5cm}

\centerline{{\bf 1 Introduction}}

\vspace{0.3cm}

In October of 2012, Adachi, Iyama and Reiten introduced the notion
of $\tau$-tilting modules which is a generalization of the classical
tilting modules [APR, BB, HR]. $\tau$-tilting modules which admit
very similar properties to the classical tilting modules are very
close to silting objects in [AiI] and the cluster tilting objects in
2-Calabi-Yau triangulated categories [IY]. So it is interesting to
find $\tau$-tilting modules for a given algebra. It is showed in
[AIR] that all $\tau$-tilting modules can be written as finite
copies of direct sums of $\tau$-rigid modules which was firstly
introduced in [AuS]. To find the $\tau$-tilting modules for given
algebras, what we need to do is just to find the (indecomposable)
$\tau$-rigid modules for them.

Notice that Adachi, Iyama and Reiten showed that every $\tau$-rigid
module $M$ has no common non-zero direct summands in the first and
second projective terms of its minimal projective resolution. It is
interesting to consider whether the $\tau$-rigid modules can be
determined by the non-existence of common direct summands in the
first and second projective terms of their minimal projective
resolutions. A positive answer to this question would make us be
able to judge $\tau$-rigid modules straightly. Unfortunately, it is
far from being true. So we have to ask: (1) When can $\tau$-rigid
modules be determined by the non-existence of common direct summands
in their minimal projective resolution?

In addition, what we also want to know is to determine the
structures of algebras from the properties of their $\tau$-rigid
modules. It is well-known that a local algebra admits a unique
$\tau$-tilting module, that is, all indecomposable $\tau$-rigid
modules are projective. So it is natural to ask: (2) Is an algebra
$\Lambda$ local if it admits a unique $\tau$-tilting module? We
should remark that a similar question for the classical tilting
modules is not true in general since every non-local self-injective
algebra admits a unique classical tilting module.

On the other hand, algebras with radical square zero have been
studied by Auslander, Reiten and Smal$\phi$ in [AuRS], which play an
important role in classifying Nakayama algebras and stable
equivalence. For the recent development of this class of algebras,
we refer to [C] and [RX]. We should note that this kind of algebras
make us be able to give more examples for algebras with best
properties of $\tau$-rigid modules and non-trivial CM-free algebras
(all finitely generated indecomposable Gorenstein projective modules
are projective).

In this paper, we try to answer the two questions above over
algebras with radical square zero. The paper is organized as
follows:

In Section 2, we will recall some preliminaries on algebras with
radical square zero. In Section 3, we give an answer to the first
question above and prove the following:

\noindent{\bf Theorem 1} {\it Let $\Lambda$ be a basic and connected
Nakayama algebra with $r^2=0$ which is not self-injective local and
let $n$ be the number of non-isomorphic simple modules. Then

(1) Every indecomposable module $M$ is $\tau$-rigid.

(2) Every $\tau$-tilting module $T$ is of the form $S_1\bigoplus
S_2\bigoplus\cdots S_t\bigoplus (\Lambda/ P_0(\tau
(S_1\bigoplus\cdots\bigoplus S_t))$, where $S_j$ is simple for
$1\leq j\leq t$, $t$ is an integer such that $0\leq t\leq {\rm
int}(n/2)$ and ${\rm int}(m)$ denotes the largest integer less than
or equal to $m$ for any real number $m$.}

\vspace{0.2 cm}

\noindent{\bf Theorem 2} {\it Let $\Lambda$ be a basic and connected
algebra with $r^2=0$.

(1) If $\Lambda$ is self-injective local, then every indecomposable
$\tau$-rigid module is projective.

(2) If $\Lambda $ is self-injective but not local, then every
indecomposable module is $\tau$-rigid.

(3) Let $M$ be an indecomposable $\Lambda$-module. If $\Lambda$ is
representation finite of finite global dimension or $\tau M$ is
$\tau$-rigid, then $M$ is $\tau$-rigid if and only if there is no
non-zero direct summand of $P_0(M)$ and $P_1(M)$ in common, where
$P_0(M)$ and $P_1(M)$ are the first and second projective terms of a
minimal projective resolution of $M$, respectively.}

\vspace{0.2cm}

In Section 4, we will give a construction theorem to get
indecomposable $\tau$-rigid modules from simple modules. This is
very different from the mutation theorem in [AIR]. As a result, we
can give an answer to the second question and prove the following:

\vspace{0.2 cm}

\noindent{\bf Theorem 3} {\it Let $\Lambda$ be a basic and connected
algebra with $r^2=0$. If $\Lambda$ admits a unique $\tau$-tilting
module, then it is local.}

\vspace{0.2cm}

In Section 5, we will give examples to show our results.

Throughout this paper, all algebras are basic connected
non-semi-simple Artin algebras over a commutative Artin ring R.
$\mathbb{D}={\rm Hom}_R(-,I^0(R/r))$ is  the ordinary dual, where
$r$ is the Jacobson radical of $R$ and $I^0(R/r)$ is the injective
envelope of $R/r$. All modules are finitely generated left
$\Lambda$-modules if not claimed.

\vspace{0.5cm}

\centerline{{\bf 2 Properties for algebras with radical square
zero}}

\vspace{0.3cm}

In this section we will recall some properties for algebras with
radical square zero. Denote by $r$ the Jacobson radical of an
algebra $\Lambda$. $\Lambda$ is called radical square zero if
$r^2=0$. Let $\Gamma$ be another algebra. We say that $\Lambda$ is
stable equivalent to $\Gamma$ if there is an equivalence functor
$F:\underline{\mod}\Lambda\rightarrow \underline{\mod}\Gamma$, where
$\underline{\mod}\Lambda$  and $\underline{\mod}\Gamma$ denote the
associate module categories modulo the projective modules,
respectively.

Now we can recall the following result for algebras with radical
square zero from [AuRS, X, Theorem 2.4, Lemma 2.1].

\vspace{0.2cm}

\noindent{\bf Lemma 2.1}  {\it Let $\Lambda$ be an algebra with
$r^2=0$. Denote by $\Gamma$ the triangular matrix algebra
\begin{equation} \left(
\begin{array}{cc}
\Lambda/r & 0 \\
r &\Lambda/r

\end{array}
\right)
\end{equation} and denote by $F: \underline{\mod}\Lambda\rightarrow
\underline{\mod}\Gamma$ the functor via $F(M)=(M/rM, rM, f)$ and
$F(g)=(g_1,g_2)$ for any $M,N,L\in\mod\Lambda$ and $g:N\rightarrow
L$, where $f:r\bigotimes_{\Lambda/r}M/rM\rightarrow rM$ is an
epimorphism, $g_1:N/rN\rightarrow L/rL$ and $g_2:rN\rightarrow rL$
are induced by $g$. Then

(1) $F$ is an equivalence and hence $\Lambda$ is stable equivalent
to $\Gamma$.

(2) $F(M)$ is indecomposable if and only if $M$ is indecomposable.

(3) $F(M)$ is projective if and only if $M$ is projective. }

\vspace{0.2cm}

Recall that a morphism $h:E\rightarrow M$ is called right minimal if
for any $l: E\rightarrow E$ $h=hl$ implies that $l$ is an
isomorphism. $h$ is right almost split if $h$ is not a spit
epimorphism and for any $m:N\rightarrow E$ which is not a split
epimorphism there exists a $t:N\rightarrow E$ such that $m=ht$.
Dually, one can define left minimal morphisms and left almost split
sequences. An exact sequence $0\rightarrow
A\stackrel{g}{\rightarrow} B\stackrel{h}{\rightarrow}C\rightarrow 0$
is called almost split if $g$ is left almost split and $h$ is right
almost split. Now we are ready to recall the following properties of
almost split sequences for algebras with $r^2=0$ from [AuRS, V,
Proposition 3.5, X, Proposition 2.5].

\vspace{0.2cm}

\noindent{\bf Lemma 2.2} {\it Let $\Lambda$ be an algebra with
$r^2$=0 and let $0\rightarrow A\stackrel{g}{\rightarrow}
B\stackrel{h}{\rightarrow}C\rightarrow 0   $ be an almost split
sequence. Then

(1) $B$ is projective if and only if $A$ is non-injective simple. If
$A$ is simple non-injective, then $h:B\rightarrow C$ is a projective
cover.

(2) $B$ is injective if and only if $C$ is non-projective simple. If
$C$ is simple non-projective, then $g:A\rightarrow B$ is an
injective envelope. }

\vspace{0.2cm}

\noindent{\bf Lemma 2.3} {\it Let $\Lambda$, $F$ and $\Gamma$ be as
in Lemma 2.1 and let $0\rightarrow A\stackrel{g}{\rightarrow}
B\stackrel{h}{\rightarrow}C\rightarrow 0\ \ (*)$ be an exact
sequence such that $A$ and $C$ are indecomposable and $A$ is not
simple. Then

(1) The sequence $(*)$ is almost split in $\mod\Lambda$ if and only
if $0\rightarrow F(A)\stackrel{F(g)}{\rightarrow}
F(B)\stackrel{F(h)}{\rightarrow}F(C)\rightarrow 0$ is almost split
in $\mod\Gamma$.

(2) If $(*)$ is almost split, then $F(A)=F(\tau_{\Lambda}
C)=\tau_\Gamma F(C)$.}

\vspace{0.2cm}

\noindent{\it Proof.} (2) follows from (1). $\hfill{\square}$

\vspace{0.2cm}

The following result which gives a connection between morphisms in
$\mod\Lambda$ and $\mod\Gamma$ is very important to the proof of the
main results.

\noindent{\bf Lemma 2.4}  {\it Let $\Lambda$, $F$ and $\Gamma$ be as
in Lemma 2.1. Then

 (1) For any $M,N\in\mod\Lambda$ we have the following
exact sequence of Abelian groups:
$$0\rightarrow {\rm Hom}_{\Lambda}(M,rN)\rightarrow {\rm Hom}_{\Lambda}(M,N)\rightarrow {\rm Hom}_{\Gamma}(F(M),F(N))\rightarrow 0$$

(2) $\underline{\rm Hom}_{\Lambda}(M,N)\simeq {\rm
Hom}_{\Gamma}(F(M),F(N))$ if both $M$ and $N$ have no projective
direct summands.}

\vspace{0.2cm}

\noindent{\it Proof.} (1) follows from [AuRS, X, Lemma 2.1] and (2)
follows from [AuRS, X, Lemma 2.3] and (1). $\hfill{\square}$

\vspace{0.2cm}

In order to show the main result on Nakayama algebra with $r^2=0$,
we need the following:

\vspace{0.2cm}

\noindent{\bf Lemma 2.5} {\it Let $\Lambda$ be a Nakayama algebra
with $r^2=0$. Then every indecomposable module $M\in \mod\Lambda$ is
either simple or projective. }

\vspace{0.2cm}

\noindent{\it Proof.} By using $r^2=0$ and [AsSS, V, Theorem 4.1].
$\hfill{\square}$

\vspace{0.5cm}

\centerline{\bf 3 $\tau$-rigid modules and minimal projective
resolution}

\vspace{0.2cm}

In this section, we will determine the $\tau$-rigid modules in terms
of minimal projective resolution and try to answer the first
question (see Theorem 3.4, Theorem 3.12 and Theorem 3.15). Firstly,
we recall the notions of $\tau$-tilting modules and $\tau$-rigid
modules in [AIR] and [AuS], respectively.

\vspace{0.2cm}

\noindent{\bf Definition 3.1} For an algebra $\Lambda$,  a
$\Lambda$-module $M$ is called $\tau$-rigid if ${\rm Hom}(M,\tau
M)=0$, where $\tau$ denotes the Auslander-Reiten translation. A
module $N$ is $\tau$-tilting if it is $\tau$-rigid and
$|N|=|\Lambda|$,  where $|N|$ denotes the number of non-isomorphic
direct summands of $N$. Any $\tau$-rigid module is a direct summand
of a $\tau$-tilting module. We also note that if $\Lambda$ is
hereditary then $\tau$-tilting modules and $\tau$-rigid modules
coincide with tilting modules and rigid modules, respectively.

For any indecomposable $M$ in $\mod\Lambda$, if $M$ is projective,
then it is $\tau$-rigid. So we can assume that $M$ is not
projective. Denote by $\cdots\rightarrow P_t(M)\rightarrow
\cdots\rightarrow P_1(M)\rightarrow P_0(M)\rightarrow M\rightarrow0$
be a minimal projective resolution of $M$, where $t$ is a
non-negative integer. And denote by $\Omega^i M$ the $i$-th syzygy
of $M$ for any $i\geq 0$. Considering the almost split sequence
$0\rightarrow \tau M\rightarrow E\rightarrow M\rightarrow 0$, we
have the following:

\vspace{0.2cm}

\noindent{\bf Proposition 3.2} {\it Let $\Lambda$ be an algebra with
$r^2=0$. Then $\Lambda$ is self-injective local if and only if there
is an almost split sequence $0\rightarrow S\rightarrow P\rightarrow
S\rightarrow 0$, where $S$ is a simple $\Lambda$-module. }

\vspace{0.2cm}

\noindent{\it Proof.} $\Rightarrow$ Since $\Lambda$ is a basic
connected local algebra, one can get the following exact sequence:
$0\rightarrow r\rightarrow\Lambda\rightarrow S\rightarrow0 ,$ where
$S$ is simple and $r$ is the radical of $\Lambda$. Notice that
$r^2=0$ then $r$ is semi-simple. Because $\Lambda$ is
self-injective, we get $r$ is simple by [HuZ, Lemma 2.6], and hence
$r\simeq S$. Then the sequence is almost split by [AsSS, IV,
Proposition 3.11].

$\Leftarrow$ By Lemma 2.2, one gets that $P$ is projective and
injective. It is enough to prove that $\Lambda$ has a unique simple
module $S$ up to isomorphism. On the contrary, Suppose that there is
another simple $S'\not\simeq S$. We claim that ${\rm
Hom}_{\Lambda}(P_0(S'),P_0(S))={\rm
Hom}_{\Lambda}(P_0(S),P_0(S'))=0$, where $P_0(M)$ is the projective
cover of $M$.

(1) ${\rm Hom}_{\Lambda}(P_0(S'),P_0(S))=0$.

Suppose that there is an $f\in {\rm Hom}_{\Lambda}(P_0(S'),P_0(S))$
such that $f\not=0$, then $f$ is not epic since $P_0(S)$ is
projective and $S\not\simeq S'$. Denote by ${\rm Im} f$ the image of
$f$, then ${\rm Im}f\subseteq rP_0(S)$. Notice that  $0\rightarrow
S\rightarrow P\rightarrow S\rightarrow 0$ is almost split, then
$P_0(S)\simeq P$ and $rP_0(S)\simeq S$ by Lemma 2.2, and hence ${\rm
Im} f=S$, then $f:P_0(S')\rightarrow S$ is epic, and hence
$P_0(S')\simeq P_0(S)$. One gets a contradiction since $S\not\simeq
S'$.

(2) ${\rm Hom}_{\Lambda}(P_0(S),P_0(S'))=0$.

Suppose that there is a $g\in {\rm Hom}_{\Lambda}(P_0(S),P_0(S'))$
such that $g\not=0$, then $g$ is not epic and ${\rm Im}g \subseteq
rP_0(S')$ by a similar argument in (1). Notice that $r^2=0$, then
$rP_0(S')$ is semi-simple. So we get
 ${\rm Im}g=S$, and hence $j: S\hookrightarrow P_0(S')$. Then we have
 the following commutative diagram:

 $$\xymatrix{\ & P_0(S')& & & \\
0\ar[r]&S{\ar[u]^j}\ar[r]^{i}&P\ar@{-->}[lu]_{h}\ar[r]&S\ar[r]&0
}.$$

By Lemma 2.2 we get that $P_0(S)\simeq P, $ $P$ is injective and
$i:S\rightarrow P$ is an injective envelope. Then $i$ is an
essential monomorphism implies that $h: P\rightarrow P_0(S')$ is
monic, and hence $P\simeq P_0(S')$, a contradiction.

Since $\Lambda$ is connected, one gets the assertion by the claim
and [AsSS, II, Lemma 1.6].  $\hfill\square$

\vspace{0.2cm}

Now we can give a class of $\tau$-rigid modules over algebras with
$r^2=0$ which are not local-self-injective.

\vspace{0.2cm}

\noindent{\bf Proposition 3.3} {\it Let $\Lambda$ be an
 algebra with $r^2=0$ which is not self-injective local and let $M$
 be an indecomposable $\Lambda$-module.
 If $M$ satisfies (1) $\tau M$ is simple projective, or (2) both $\tau M$ and $M$ are
 simple, then $M$ is $\tau$-rigid.}

 \vspace{0.2cm}

\noindent{\it Proof.} Suppose that there is a non-zero $f\in$ ${\rm
Hom}_{\Lambda}(M,\tau M)$. For both cases, we have $f$ is epic, and
hence $M\simeq\tau M$. For the first case, we get that $M$ is
projective, a contradiction. For the second one, by Lemma 2.2 there
is an almost split sequence $0\rightarrow S\rightarrow P\rightarrow
S\rightarrow 0$ with $S$ a simple module. By Proposition 3.2
$\Lambda$ is self-injective local, a contradiction. $\hfill \square$

\vspace{0.2cm}

\noindent{\bf Remark} For any algebra $\Sigma$ and an indecomposable
$\Sigma$-module $N$ with $\tau N$ simple projective, one can show
that $N$ is $\tau$-rigid by formulating the proof of Proposition 3.3
(1).

Now we are in a position to state the $\tau$-tilting and
$\tau$-rigid modules for Nakayama algebras with $r^2=0$.

\vspace{0.2cm}

\noindent{\bf Theorem 3.4} {\it Let $\Lambda$ be a Nakayama algebra
with $r^2=0$ which is not self-injective local and let $n$ be the
number of non-isomorphic simple modules. Then

(1) Every indecomposable module $M$ is $\tau$-rigid.

(2) Every $\tau$-tilting module $T$ is of the form $S_1\bigoplus
S_2\bigoplus\cdots S_t\bigoplus (\Lambda/ P_0(\tau
(S_1\bigoplus\cdots\bigoplus S_t))$, where $S_j$ is simple for
$1\leq j\leq t$, $t$ is an integer such that $0\leq t\leq {\rm
int}(n/2)$ and ${\rm int}(m)$ denotes the largest integer less than
or equal to $m$ for any real number $m$.}

\vspace{0.2cm}

\noindent{\it Proof.} (1) By Lemma 2.5, $M$ is simple or projective
for any indecomposable $M\in\mod\Lambda$. If $M$ is projective,
there is nothing to prove. If $M$ is simple non-projective, then
$\tau M$ is simple, by Proposition 3.3, $M$ is $\tau$-rigid.

(2) For any $\tau$-tilting module $T$, we claim that there is at
least one indecomposable $P$ as a direct summand of $T$.

On the contrary, suppose that $T=S_1\bigoplus
S_2\bigoplus\cdots\bigoplus S_n$ with all $S_j$ simple
non-projective for $1\leq j\leq n$. Without loss of generality, we
can assume that $\tau S_1=S_2$ by Lemma 2.5 or (1). Then $0\not={\rm
Hom}_{\Lambda}(S_2,S_2)\subseteq{\rm Hom}_{\Lambda}(T,\tau T)=0$, a
contradiction.

Next we will show if $S$ is a direct summand of $T$, then $P_0(S)$
is a direct summand of $T$.

If $S$ is projective, then the assertion holds true. We can assume
that $S$ is not projective. Then we get the following almost split
sequence: $0\rightarrow S'\rightarrow P\rightarrow S\rightarrow 0$
with $S'$ simple by Lemma 2.2 and Lemma 2.5. Again by Lemma 2.2, $P$
is indecomposable projective and $P\simeq P_0(S)$. Because $T$ is
$\tau$-tilting and $S$ is direct summand of $T$, it is not difficult
to show that $S'$ is not a direct summand of $T$. Similarly, if $S$
is not injective, then a simple module $S^*$ with $\tau S^*\simeq S$
is not a direct summand of $T$. Since $\tau T$ is semi-simple, then
${\rm Hom}_{\Lambda}(P_0(S),\tau T)=0$, that is, $P_0(S)$ is in
$Fac(T)$  by [AIR, Theorem 2.10], where ${\rm Fac(T)}$ denotes the
category consisting of factor modules of finite copies of direct
sums of $T$. So $P_0(S)$ is a direct summand of $T$. Notice that
$|T|=|\Lambda|=n$, so the number of simple direct summands of $T$
has to be at most ${\rm int(n/2)}$. $\hfill\square$

\vspace{0.2cm}

The conditions in Proposition 3.3 are not easy to be satisfied. In
the following we will generalize it into a general framework. Denote
by $\gldim\Lambda$ the global dimension of $\Lambda$ and denote by
$\pd_{\Lambda}M$ the projective dimension of $M$. We have:

\vspace{0.2cm}

\noindent{\bf Lemma 3.5} {\it Let $\Lambda$ be an algebra with
$r^2=0$. If $S$ is a simple module with $\pd_{\Lambda}S=m<\infty$.
Then $S$ is $\tau$-rigid. Moreover, if $\gldim\Lambda=m<\infty$,
then every simple module $S$ is $\tau$-rigid.}

\vspace{0.2cm}

\noindent{\it Proof.} We only have to the first one since the last
follows from the first. By [AIR, Proposition 1.2] a simple module
$S$ is $\tau$-rigid if and only if it is rigid, that is, ${\rm
Ext}_{\Lambda}^1(S,S)=0$. If $S$ is projective, there is nothing to
prove. So we can assume that $m\geq1$. Take the following part of a
minimal projective resolution of $S$: $0\rightarrow \Omega^1
S\rightarrow P_0(S)\rightarrow S\rightarrow 0$, where $\Omega^1 S$
denotes the first syzygy of $S$. One gets that
$\pd_{\Lambda}\Omega^1 S=m-1$ since $\pd_{\Lambda}S=m<\infty$. Since
$r^2=0$, we have $\Omega^1 S\simeq rP_0(S)$ is semi-simple and any
direct summand of it is of projective dimension at most $m-1$. So it
is not difficult to show ${\rm Ext}_{\Lambda}^{1}(S,S)\simeq {\rm
Hom}_{\Lambda}(\Omega^1 S,S) =0.$ $\hfill{\square}$

\vspace{0.2cm}

Denote by $\overline{\mod}\Lambda$ the associate modules category
modulo injective modules and denote by $\overline{{\rm
Hom}}_\Lambda(L,N)$ and $\underline{{\rm Hom}}_\Lambda(L,N)$ classes
of morphisms from $L$ to $N$ in $\overline{\mod}\Lambda$ and
$\underline{\mod}\Lambda$, respectively. Now we are in a position to
state another main result on judging the $\tau$-rigid properties by
simple modules.

\vspace{0.2cm}

\noindent{\bf Theorem 3.6} {\it  Let $\Lambda$ be an algebra with
$r^2=0$ and let $M$ be indecomposable with $\tau M$ simple. We have
(1) $M$ is $\tau$-rigid if and only if $\tau M$ is $\tau$-rigid. (2)
If $\pd_{\Lambda}M<\infty$, then $M$ is $\tau$-rigid. Moreover, if
$\gldim\Lambda<\infty$, then $M$ is $\tau$-rigid. }

\vspace{0.2cm}

\noindent{\it Proof.} Since (2) is a straight result of (1), Lemma
2.2 and Lemma 3.5, we only show (1).

$\Leftarrow$ By the remark of Proposition 3.3, it is enough to show
the case of $\tau M$ is not projective.

On the contrary, suppose that $M$ is not $\tau$-rigid, that is,
${\rm Hom}(M,\tau M)\not=0$. We get that $f$ is epic for any
$0\not=f\in{\rm Hom}(M,\tau M)$ since $\tau M$ is simple. By Lemma
2.2, one gets the following almost split sequence: $0\rightarrow\tau
M\stackrel{i}{\rightarrow} P_0(M)\rightarrow M\rightarrow0$. So
$P_0(\tau M)$ is a direct summand of $P_0(M)$. Notice that $i$ is
left minimal, then by [AuRS, I, Theorem 2.4] $\tau M$ can be
embedded into $P_0(\tau M)$, and hence a direct summand of $r
P_0(\tau M)$ since $r^2=0$. Then we have the following commutative
diagram:

$$\xymatrix{& 0 \ar[r] & rP_0(\tau M)\ar[r]
\ar[d]^{\alpha} & P_{0}(\tau M) \ar[r]^{\pi}\ar[d]^{\beta} & \tau M \ar[r]\ar@{=}[d] & 0\\
& 0 \ar[r] & \tau M\ar[r] & E \ar[r] & \tau M \ar[r] & 0}$$
 with $\alpha$ an epimorphism. By the snake lemma, $\beta$ is also
 an epimorphism. Since $\pd_{\Lambda}M<\infty$, one can show that $\pd_{\Lambda}\tau
 M$ is of finite projective dimension by Lemma 2.2. Then by the assumption $\tau M$ is
 $\tau$-rigid and hence the bottom row in the commutative diagram is
 split. So one gets an epimorphism $P_0(\tau M)\rightarrow \tau M\bigoplus \tau
 M$, a contradiction.

 $\Rightarrow$ Since $M$ is $\tau$-rigid, one gets ${\rm Hom}_{\Lambda}(M,\tau
 M)=0$ which implies that $\underline{{\rm Hom}}_{\Lambda}(M,\tau
 M)=0$. Notice that $\tau:
 \underline\mod\Lambda\rightarrow\overline\mod\Lambda$ is an
 equivalence, one can get $\overline{{\rm Hom}}_{\Lambda}(\tau
 M,\tau^2 M)=\underline{{\rm Hom}}_{\Lambda}(M,\tau
 M)=0$. By AR-formula one gets ${\rm Ext}_{\Lambda}^1(\tau M,\tau M)\simeq \mathbb{D}\overline{{\rm Hom}}_{\Lambda}(\tau M,\tau^2 M)=0$.
  Then by [AIR, Proposition 1.2] $\tau M$ is $\tau$-rigid since
$\tau M$ is simple. $\hfill\square$

\vspace{0.2cm}

To answer the first question, we need the following properties for
$\tau$-rigid modules over hereditary algebras.

\vspace{0.2cm}

\noindent{\bf Lemma 3.7} {\it Let $\Lambda$ be a hereditary algebra
and let $M$ be an indecomposable non-projective module. If $\tau M$
is projective, then $M$ is $\tau$-rigid. }

\vspace{0.2cm}

\noindent{\it Proof.} Suppose that ${\rm Hom}_{\Lambda}(M,\tau
M)\not=0$. Then there is a non-zero morphism $f:M\rightarrow \tau
M$. Since $\Lambda$ is hereditary and $\tau M$ is projective, one
gets that ${\rm Im} f$ is projective and hence $M$ is projective, a
contradiction. $\hfill\square$

\vspace{0.2cm}

\noindent{\bf Lemma 3.8} {\it Let $\Lambda$ be a hereditary algebra
and let $M$ be an indecomposable non-projective module. If $\tau M$
is $\tau$-rigid, then $M$ is $\tau$-rigid.}

\vspace{0.2cm}

\noindent{\it Proof.} If $\tau M$ is projective, the assertion holds
from Lemma 3.7. We only show the case $\tau M$ is not projective. By
[AsSS, IV, Corollary 2.15 (b)], one gets ${\rm Hom}_{\Lambda}(M,\tau
M)\simeq {\rm Hom}_{\Lambda}(\tau M, \tau^2 M)=0$.   $\hfill
\square$

\vspace{0.2cm}

Recall that an indecomposable module $M$ over a hereditary algebra
is preprojective if there is a non-negative integer $j$ such that
$\tau^j M$ is a non-zero projective module. Then we have:

 \vspace{0.2cm}

\noindent{\bf Proposition 3.9} [AuRS, VIII, Propositions 1.7, 1.13]
{\it Let $\Lambda$ be a hereditary algebra. Then (1) Every
preprojective module $M$ is $\tau$-rigid. (2) If $\Lambda$ is of
finite type, then every indecomposable module is $\tau$-rigid, and
hence rigid.}

\vspace{0.2cm}

\noindent{\it Proof.} (1) We can assume that $\tau^j M$ is
projective for some non-negative integer $j$. By induction on $j$
and Lemma 3.8, one gets the assertion. Then by [AuRS, VIII,
Proposition 1.13] and (1), one can show (2). $\hfill\square$

\vspace{0.2cm}

 Denote by $\id_{\Lambda}M$ $(resp.\ \id_{\Lambda^o}M)$ the
injective dimension of $M$ for an $M$ in $\mod\Lambda$ $(resp.\
\mod\Lambda^o)$. Recall that an algebra $\Lambda$ is called
Gorenstein if $\id_{\Lambda}\Lambda=\id_{\Lambda^o}\Lambda=n$ for
some integer $n\geq 0$. We have the following:

\vspace{0.2cm}

\noindent{\bf Lemma 3.10} {\it Let $\Lambda$ be a Gorenstein algebra
with $r^2=0$. Then $\Lambda$ is either self-injective or of finite
global dimension.}

\vspace{0.2cm}

\noindent{\it Proof.} By [C], we can get that every algebra with
$r^2=0$ is either self-injective or CM-free. Recall that an algebra
is called CM-free if every finitely generated Gorenstein projective
module is projective. We will show $\gldim\Lambda=n$ if
$\id_{\Lambda}\Lambda=\id_{\Lambda^o}\Lambda=n$ for some $n>0$. By
[AuR, Proposition 3.1] one can show that $\Omega^n M$ is Gorenstein
projective and hence projective for any $M$ in $\mod\Lambda$ since
$\Lambda$ is CM-free. The assertion holds true. $\hfill\square$

\vspace{0.2cm}

Notice that a self-injective algebra with $r^2=0$ is Nakayama by
[AuRS, IV, Proposition 2.16]. Since the Nakayama case is completely
classified in Proposition 3.2 and Theorem 3.4, we only have to find
the $\tau$-rigid modules for algebras of finite global dimension
with $r^2=0$.

 Denote by ${\rm Soc} M$
the socle of $M$. We have the following easy observation:

\vspace{0.2cm}

\noindent{\bf Lemma 3.11} {\it Let $\Lambda$ be an algebra and let
$M$ be an indecomposable and non-projective $\Lambda$-module. Then
(1) ${\rm Soc} \tau M\simeq \Omega^1 M/r\Omega^1 M$. (2) ${\rm Soc}
\tau M\simeq \Omega^1 M$ if $r^2=0$.}

\vspace{0.2cm}

\noindent{\it Proof.} (1) Taking the following part of a minimal
projective resolution of $M$: $P_1(M)\rightarrow P_0(M)\rightarrow
M\rightarrow0$ and then applying the functor ${\rm
Hom}_\Lambda(-,\Lambda)=(-)^*$, one gets the following part of
minimal projective resolution
 of ${\rm Tr} M$: $ P_{0}(M)^*\rightarrow P_{1}(M)^*\rightarrow {\rm Tr} M\rightarrow 0$, where ${\rm Tr}$ is the Auslander-Bridger transpose.
Then applying the functor $\mathbb{D}$, one has a minimal injective
resolution of $\tau M:$ $0\rightarrow\tau M\rightarrow
\mathbb{D}P_1(M)^*\rightarrow \mathbb{D}P_0(M)^*$.  Then ${\rm
Soc}\tau M\simeq {\rm Soc}\mathbb{D}P_1(M)^*\simeq \Omega^1
M/r\Omega^1 M.$

(2) Taking the following part of a minimal projective resolution of
$M$: $0\rightarrow \Omega^1 M\rightarrow P_0(M)\rightarrow
M\rightarrow0$, one can get the following commutative diagram:
$${\xymatrix{0\ar[r]&\Omega^1 M\ar[r]\ar[d]^g&P_0(M)\ar[r]\ar@{=}[d]&M\ar[r]\ar[d]^f&0\\
0\ar[r]&rP_0(M)\ar[r]&P_0(M)\ar[r]&M/rM\ar[r]&0}}, $$ where $f$ is
epic. Then by the snake lemma, one gets $g$ is a monomorphism. Since
$r^2=0$, we have that $rP_0(M)$ is semi-simple and hence $\Omega^1
M$ is semi-simple. We are done. $\hfill\square$

\vspace{0.2cm}

Now we are ready to determine all the $\tau$-rigid modules for a
Gorenstein algebra of finite type with $r^2=0$. Combined with
Proposition 3.2, Lemma 3.11 and Theorem 3.4, one can show:

\vspace{0.2cm}

\noindent{\bf Theorem 3.12} {\it Let $\Lambda$ be a Gorenstein
algebra of finite type with $r^2=0$.

(1) If $\Lambda $ is self-injective local, then every indecomposable
$\tau$-rigid module is projective.

(2) If $\Lambda $ is self-injective but not local, then every
indecomposable module $M$ is $\tau$-rigid

(3) If $\Lambda$ is of finite global dimension, then an
indecomposable module $M$ is $\tau$-rigid if and only if there is no
non-zero direct summand of $P_0(M)$ and $P_1(M)$ in common.}

\vspace{0.2cm}

\noindent{\it Proof.} (1) is clear and (2) is showed in Theorem 3.4.

(3) $\Rightarrow$ It is a straight result of [AIR, Proposition 2.5].

$\Leftarrow$ Without loss of generality, we can assume that $\tau M$
is not zero. If $\tau M$ is simple, then the assertion holds true by
Theorem 3.6. Now we can assume that $\tau M$ is not simple. Let
$\Gamma$ and $F$ be as in Lemma 2.1. By Lemma 2.1 and [AuRS, X,
Proposition 1.1], we get that $\Gamma$ is hereditary of finite type.
Then by Proposition 3.9(2), $F(M)\in\mod\Gamma$ is $\tau$-rigid. So
${\rm Hom}_{\Gamma} (F(M),F(\tau M))\simeq {\rm Hom}_{\Gamma}
(F(M),\tau F(M))=0$ by Lemma 2.3. Then by Lemma 2.4 $M$ is
$\tau$-rigid if and only if ${\rm Hom}_{\Lambda}(M,r\tau M)=0$.

We show that ${\rm Hom}_{\Lambda}(M, r\tau M)=0$. Since $r^2=0$, we
get that $r\tau M$ is semi-simple and hence a direct summand of
${\rm Soc} \tau M\simeq \Omega^1 M$ by Lemma 3.11. Notice that there
is no common direct summand of $P_0(M)$ and $P_1(M)$, one can show
${\rm Hom}_{\Lambda}(M, {\rm Soc}\tau M)=0$ which implies that ${\rm
Hom}_{\Lambda}(M, r\tau M)=0$. Then $M$ is $\tau$-rigid by Lemma
2.4. $\hfill\square$

\vspace{0.2 cm}

In general, for algebras mentioned in Theorem 3.12 (3) we don't know
whether there is a common direct summand in $P_0(M)$ and $P_1(M)$
for an indecomposable $M$ (see Example 5.3). However, we get the
following:

\vspace{0.2cm}

\noindent{\bf Proposition 3.13} {\it  Let $\Lambda$ be an algebra of
finite global dimension with $r^2=0$ and let $M$ be an
indecomposable module.
 If $M/rM$ is simple, then $P_0(M)$ and $P_1(M)$ have no non-zero
 direct summands in common.}

\vspace{0.2cm}

\noindent{\it Proof.} Denote by $S=M/rM$. It is enough to show that
${\rm Hom}_{\Lambda}(\Omega^1 M, S)=0$. On the contrary, suppose
that ${\rm Hom}_{\Lambda}(\Omega^1 M, S)\not=0$. Then $S$ is a
direct summand of $\Omega^1 M$ since $r^2=0$ and $S$ is simple.
Moreover, we have the following commutative diagram with exact rows:

$$\xymatrix{0\ar[r]& \Omega^1 M\ar[r]\ar[d]^f& P_0(M)\ar[r]\ar@{=}[d] & M\ar[r]\ar[d]^g&0 \\
0\ar[r]&\Omega^1 S\ar[r]&P_0(S)\ar[r]&S\ar[r]&0 }.$$ Since $g$ is
epic, one can show that $f:\Omega^1 M\rightarrow \Omega^1 S$ is a
monomorphism, and hence $\Omega^1 M$ is a direct summand of
$\Omega^1 S$ because $r^2=0$. So $S$ is a direct summand of
$\Omega^1 S$, that is $\pd_{\Lambda}S\leq \pd_{\Lambda}S-1<\infty$,
a contradiction. $\hfill\square$

\vspace{0.2cm}

Denote by $l(M)$ the length of $M$. As a result of Proposition 3.13,
we can get:

\vspace{0.2cm}

\noindent{\bf Corollary 3.14} {\it Let $\Lambda$ be an algebra of
finite global dimension with $r^2=0$ and $M$ be an indecomposable
module. if $l(M)\leq 2$, then $P_0(M)$ and $P_1(M)$ have no non-zero
 direct summands in common. }

\vspace{0.2cm}

\noindent{\it Proof.} By Proposition 3.13, it suffice to prove that
$M/rM$ is simple. If $l(M)=1$ then $M$ is simple, the assertion
holds true. If $l(M)=2$ then $rM={\rm Soc} M$ is simple, and hence
$M/rM$ is simple. $\hfill\square$

\vspace{0.2cm}

At the end of this section we will give a method to find more
$\tau$-rigid modules for any algebra with $r^2=0$. As we know all
the projective $\Lambda$-modules are $\tau$-rigid. In case that the
algebra $\Lambda$ is not self-injective, there must be some
indecomposable $M$ such that $\tau M$ is projective. It is
interesting to know whether $M$ is $\tau$-rigid. A more general
question is: Whether is $M$ $\tau$-rigid if $\tau M$ is
 $\tau$-rigid? To answer this question, we have

\vspace{0.2cm}

\noindent{\bf Theorem 3.15} {\it Let $\Lambda$ be an algebra with
$r^2=0$ and let $M$ be an indecomposable $\Lambda$-module such that
$\tau M$ is $\tau$-rigid. Then $M$ is $\tau$-rigid if and only if
there is no non-zero direct summand of $P_0(M)$ and $P_1(M)$ in
common.}

\vspace{0.2cm}

\noindent{\it Proof.} $\Rightarrow$  By [AIR, Proposition 2.5].

$\Leftarrow$  If $\tau M$ is simple, then by Theorem 3.6 (1) $M$ is
$\tau$-rigid if and only if $\tau M$ is $\tau$-rigid. Then the
assertion holds by the assumption $\tau M$ is $\tau$-rigid.

Now we can assume that $\tau M$ is not simple. Let $\Gamma$ and $F$
be as in Lemma 2.1. We claim that ${\rm Hom}_{\Gamma}(F(M),F(\tau
M))=0$.

  Since $\tau M$ is
$\tau$-rigid, we get ${\rm Hom}_{\Lambda}(\tau M,\tau^2 M)=0$ and
hence $\overline{{\rm Hom}}_{\Lambda}(\tau M,\tau^2 M)=0$. Notice
that $\tau: \underline{\mod\Lambda}\rightarrow
\overline{\mod\Lambda}$ is an equivalence, we get that
$\underline{{\rm Hom}}_{\Lambda}(M,\tau M)=0$. If $\tau M$ is not
projective, then we get ${\rm Hom}_{\Gamma}(F(M),F(\tau M))=0$ by
Lemma 2.4. If $\tau M$ is projective, then ${\rm
Hom}_{\Gamma}(F(M),F(\tau M))=0$ since $F(\tau M)$ is projective and
$\Gamma$ is hereditary. Otherwise, one can get a non-zero
$g:F(M)\rightarrow F(\tau M)$. So ${\rm Im}g$ is projective, and
hence $F(M)$ is projective, that is, $M$ is projective by Lemma 2.1,
a contradiction.

By Lemma 2.4, we only have to show ${\rm Hom}_{\Lambda}(M,r\tau
M)=0$. One can get the assertion by a similar argument in Theorem
3.12. $\hfill\square$

\vspace{0.2cm}

For a non-Nakayama algebra $\Lambda$ with $r^2=0$, by Theorem 3.15
one can find the $\tau$-rigid modules one by one from the projective
vertices of the $AR$-quiver of $\Lambda$ since here $\Lambda$ is not
self-injective. For the Nakayama case with $r^2=0$, we refer to
Theorem 3.4 or Theorem 3.12. On the other hand, it is interesting to
study the structure of algebras in terms of indecomposable
$\tau$-rigid modules. Compared with Theorem 3.10 and Lemma 2.1, we
end this section with an open question which is closed to algebras
of finite type.

{\bf Question} Let $\Lambda$ be an algebra with radical square zero.
If all the indecomposable modules are $\tau$-rigid, then $\Lambda$
is of finite type.

\vspace{0.5cm}

\centerline {\bf 4 $\tau$-rigid modules and local algebras}

\vspace{0.2cm}

 In this section we firstly introduce a theorem to get a class of indecomposable $\tau$-rigid modules from simple modules (here we don't need
 $\Lambda$ to be radical square zero). This method is very different
 from the mutation theorem in [AIR]. As a result, we give a partial
 answer to the second question.

 \vspace{0.2cm}

\noindent {\bf Theorem 4.1} {\it Let $\Lambda$ be an algebra and let
$S$ be a simple $\Lambda$-module such that the first syzygy
$\Omega^1 S$ is non-zero semi-simple.

(1) Suppose that $S$ is not a direct summand of $\Omega^1 S$. Let
$S_1$ be a simple submodule of $\Omega^1 S$ and let $m$ be the
maximal integer such that ${S_1}^m$ is a direct summand of $\Omega^1
S$. Then there is an exact sequence $0\rightarrow
\Omega^1S/S_1^m\rightarrow P_0(S)\rightarrow M\rightarrow0$ with $M$
indecomposable $\tau$-rigid.

(2) Assume that $S$ is a direct summand of $\Omega^1 S$. Let $n$ be
the maximal integer such that ${S}^n$ is a direct summand of
$\Omega^1 S$.

(a) If $\Omega^1 S\simeq S^n$, then there is no non-projective
indecomposable $\tau$-rigid module $N$ with the projective cover
$P_0(N)\simeq P_0(S)$.

(b) If $\Omega^1 S\not\simeq S^n$, then we can get the following
exact sequence $0\rightarrow \Omega^1S/S^n\rightarrow
P_0(S)\rightarrow N\rightarrow 0$ such that $N$ is indecomposable
non-projective $\tau$-rigid. }

\vspace{0.2cm}

\noindent{\it Proof.} (1) By the assumption of (1), one can get that
the simple module $S$ is non-projective $\tau$-rigid since ${\rm
Ext}_{\Lambda}^1(S,S)\simeq {\rm Hom(\Omega^1S,S)}=0$.

If $\Omega^1 S\simeq {S_1}^m$, there is nothing to prove.

Now we can assume that $\Omega^1 S\not\simeq {S_1}^m$. Since
$\Omega^1S$ is semi-simple, we get that a monomorphism
$\Omega^1S/{S_1}^m\hookrightarrow\Omega^1S\hookrightarrow P_0(S)$,
and hence we have the desired exact sequence $0\rightarrow
\Omega^1S/S_1^m\rightarrow P_0(S)\rightarrow M\rightarrow0$. It
remains to prove that $M$ is indecomposable $\tau$-rigid.

Since $P_0(S)$ is indecomposable and projective, one can show that
$M$ is indecomposable and $P_0(M)\simeq P_0(S)$ by the exact
sequence above. In the following we show that $M$ is $\tau$-rigid.
By [AIR, Proposition 1.2(a)], it is enough to show that ${\rm
Ext}_{\Lambda}^1(M,N)=0$ for any $N\in{\rm Fac}M$, where ${\rm
Fac}M$ is the full subcategory consisting of factor modules of
finite copies of direct sums of $M$.

By the construction of $M$, we have the following commutative
diagram with exact rows

$$\xymatrix{0\ar[r]& \Omega^1 M\ar[r]\ar[d]& P_0(M)\ar[r]^a\ar@{=}[d] & M\ar[r]\ar[d]&0 \\
0\ar[r]&\Omega^1 S\ar[r]&P_0(S)\ar[r]&S\ar[r]&0 }.$$ Here
$\Omega^1M\simeq \Omega^1S/S_1^m$ and by snake lemma one gets an
exact sequence $$0\rightarrow
\Omega^1M\rightarrow\Omega^1S\rightarrow S_1^m\rightarrow0\ \ \ \ \
(*1)$$

Since $N$ is in ${\rm Fac}M$, then there is a minimal positive
integer $t\geq1$ such that $g:M^t\rightarrow N$ is an epimorphism.
By [AuRS, I, Theorem 2.2] it is not difficult to show that
$P_0(N)\simeq P_0(M)^t\simeq P_0(S)^t$. Hence we get an epimorphism
$h:N\rightarrow S^t$. Then we have the following commutative diagram
with exact rows

$$\xymatrix{0\ar[r]& \Omega^1 N\ar[r]\ar[d]& P_0(N)(\simeq P_0(S)^t)\ar[r]\ar@{=}[d] & N\ar[r]\ar[d]^h&0 \\
0\ar[r]&{\Omega^1 S}^t\ar[r]&P_0(S)^t\ar[r]&S^t\ar[r]&0 }.$$ Notice
that $h$ is an epimorphism and $\Omega^1S$ is semi-simple, by snake
lemma we have two exact sequences $$0\rightarrow
\Omega^1N\rightarrow \Omega^1S^t\rightarrow L\rightarrow 0\ \ \ \
(*2)$$

$$0\rightarrow L\rightarrow N\rightarrow S^t\rightarrow 0\ \ \ (*3)$$

 On the other hand, we have the following commutative diagram

$$\xymatrix{0\ar[r]& \Omega^1 M^t\ar[r]\ar[d]^{\exists l}& P_0(M)^t(\simeq P_0(S)^t)\ar[r]\ar@{=}[d] & M^t\ar[r]\ar[d]^g&0 \\
0\ar[r]&{\Omega^1 N}\ar[r]&P_0(N)(\simeq P_0(S)^t)\ar[r]&N\ar[r]&0
}.$$ Since $g$ is an epimorphism, we get a monomorphism
$l:\Omega^1M^t\rightarrow\Omega^1N$. Combining the exact sequence
$(*1)$ and $(*2)$, we have the following commutative diagram:

$$\xymatrix{0\ar[r]& \Omega^1 M^t\ar[r]\ar[d]^l& \Omega^1S^t\ar[r]\ar@{=}[d] & S_1^{mt}\ar[r]\ar[d]^{\exists f}&0 \\
0\ar[r]&{\Omega^1 N}\ar[r]&\Omega^1S^t\ar[r]&L\ar[r]&0}.$$ By snake
lemma again, we get that $f$ is an epimorphism. Notice that
$S_1^{mt}$ is semi-simple, then $L$ is a direct summand of
$S_1^{mt}$. Applying the functor ${\rm Hom}_{\Lambda}(M,-)$ to the
exact sequence $(*3)$, one can get that ${\rm
Ext}^1_{\Lambda}(M,S)\simeq {\rm Hom}_{\Lambda}(\Omega^1M, S)=0$.
Similarly, one can get ${\rm Ext}^1_{\Lambda}(M,S_1)=0$ and hence
${\rm Ext}^1_{\Lambda}(M,L)=0$. Then one gets ${\rm
Ext}^1_{\Lambda}(M,N)=0$. We are done.

(2) We only prove (a) since the proof of (b) is very similar to the
proof of (1). It is easy to show that $S$ is not $\tau$-rigid.
Suppose that there is an indecomposable $\tau$-rigid module $N$ such
that $P_0(S)\simeq P_0(N)$. Then $N\not\simeq S$ and we have the
following commutative diagram

 $$\xymatrix{0\ar[r]& \Omega^1 N\ar[r]\ar[d]& P_0(N)\ar[r]\ar@{=}[d] & N\ar[r]\ar[d]&0 \\
0\ar[r]&\Omega^1 S\ar[r]&P_0(S)\ar[r]&S\ar[r]&0 }.$$ By snake lemma,
one get that $\Omega^1N$ is a direct summand of $\Omega^1S$, and
hence has $S$ as one of its direct summand. That means $P_0(N)$ and
$P_1(N)$ have a non-zero direct summand $P_0(S)$. But $N$ is
$\tau$-rigid, by using [AIR, Proposition 2.5], one gets a
contradiction.  $\hfill{\square}$

\vspace{0.2cm}

\noindent{\bf Remark} One can easily show that algebras with radical
square zero satisfy the condition of Theorem 4.1. For a non-local
algebra $\Gamma$ with radical square zero, there is at least $2n-m$
indecomposable $\tau$-rigid modules, where $n$ and $m$ is the number
of non-isomorphic simple modules and the number of non-isomorphic
simple projective modules, respectively.

\vspace{0.2cm}

In the following we will focus on the structure of algebras and the
homological properties of algebras for which all $\tau$-rigid
modules are projective. To prove the main result of this section, we
need the following lemmas.

\vspace{0.2cm}

\noindent{\bf Lemma 4.2 } {\it Let $\Lambda$ be an algebra such that
all $\tau$-rigid modules are projective, then $\Lambda$ has no
simple projective module.}

\noindent{\it Proof.} Suppose that there is a simple projective
module $S$. Then one can get an AR-sequence $0\rightarrow
S\rightarrow E\rightarrow M\rightarrow0$. By Proposition 3.3, $M$ is
$\tau$-rigid. But $M$ is not projective since $\tau M\simeq
S\not=0$. $\hfill{\square}$

{\vspace{0.2cm}

\noindent{\bf Lemma 4.3} {\it Let $\Lambda$ be an algebra such that
all $\tau$-rigid modules are projective and let $S$ be a simple
$\Lambda$-module.

(1)Then there is a non-zero direct summand of $P_0(S)$ and $P_1(S)$
in common.

(2) If in addition $\Lambda$ is radical square zero, then $S$ is a
direct summand of $\Omega^1S$. }

{\vspace{0.2cm}

\noindent{\it Proof.} (1) By Lemma 4.2, we get that there is no
projective simple module. By the assumption, $S$ is not
$\tau$-rigid. By [AIR, Proposition 1.2 (a)], $0\not={\rm
Ext}_{\Lambda}^{1}(S,S)\simeq{\rm Hom}_{\Lambda}(\Omega^1S,S)$. Then
one gets the assertion.

(2) is a straight result of (1). $\hfill{\square}$

{\vspace{0.2cm}

Now we can state the main theorem of this section.

\noindent{\bf Theorem 4.4} {\it Let $\Lambda$ be an algebra with
radical square zero. If $\Lambda$ admits a unique $\tau$-tilting
module, then $\Lambda$ is local. }

\noindent{\it Proof.} Firstly, we claim that for any simple module
$S$, $\Omega^1S\simeq S^t$ for some positive integer $t$. By Lemma
4.3, we get that $S$ is a direct summand of $\Omega^1S$. By Theorem
4.1 (2)(b), $\Omega^1S/S^t$ must be zero (otherwise, there will be
an indecomposable non-projective $\tau$-rigid module). The assertion
holds.

Next we will show that there is a unique simple module $S$ in
$\mod\Lambda$. Suppose there is another simple module $S'$. Then by
the claim above we get that there is a positive integer $m$ such
that $\Omega^1S'\simeq {S'}^m$. So one can get ${\rm
Hom}_{\Lambda}(P_0(S),P_0(S'))={\rm
Hom}_{\Lambda}(P_0(S'),P_0(S))=0$. Notice that $\Lambda$ is basic
and connected, this is a contradiction. $\hfill{\square}$

\vspace{0.2cm}

In Theorem 4.4, if $\Lambda$ is a finite dimensional algebra over an
algebraically closed field $K$, one can get that the quiver of
$\Lambda$ is just one vertex with several cycles. Then one
determines the structure of the algebras completely. After finishing
Theorem 4.4, the author was told by Professor Iyama that he can
prove that a basic connected algebra with a unique $\tau$-tilting
module is a local algebra by mutation.

\vspace{0.5cm}

\centerline {\bf 5 Examples}

\vspace{0.2cm}

In this section we give examples to show our results. Let
$\textbf{Q}$ be a quiver. Denote by $P(i)$, $I(i)$ and $S(i)$ the
indecomposable projective module, indecomposable injective module
and the simple module according to the vertex $i\in \textbf{Q}$,
respectively. The following example is a Nakayama algebra with
$r^2=0$.

\vspace{0.2cm}

\noindent{\bf Example 5.1} Let $\Lambda$ be given by the quiver:
$$
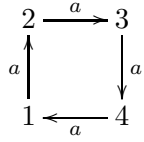
\xymatrix{2\ar[r]^a&3\ar[d]^a\\
1\ar[u]^a&4\ar[l]^a}$$ with relations $a^2=0$. By Theorem 3.4, every
indecomposable module is $\tau$-rigid. The $\tau$-tilting modules
are of the following forms:

(1) $0$-simple module.
 $P(1)\bigoplus P(2)\bigoplus P(3)\bigoplus P(4)$

 (2) $1$-simple module.
$P(1)\bigoplus S(1)\bigoplus P(3)\bigoplus P(4), P(1)\bigoplus
P(2)\bigoplus S(2)\bigoplus P(4)$ \\
\indent \ \ \ \ \ \ \ \ \ \ \ \ \ \ \ \ \ \ \ \ \ \ \ \ \ \ \ \
$P(1)\bigoplus P(2)\bigoplus P(3)\bigoplus S(3), S(4)\bigoplus
P(2)\bigoplus P(3)\bigoplus P(4)$

(3) $2$-simple modules. $P(1)\bigoplus S(1)\bigoplus P(3)\bigoplus
S(3), S(4)\bigoplus P(2)\bigoplus S(2)\bigoplus P(4)$

\vspace{0.2cm}

In the following we give an example to show that there does exist an
algebra of finite global dimension with $r^2=0$ which is of finite
type but not Nakayama.

\vspace{0.2cm}

\noindent{\bf Example 5.2} Let $\Lambda$ is given by the quiver:
 $$
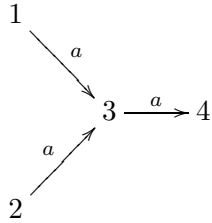
{\xymatrix{1\ar[rd]^a&\ &\ \\
 \ &3\ar[r]^a&4\\
 2\ar[ur]^a&\ &\ }}$$ with relations $a^2=0$. Then

 (1) $\Lambda$ is a representation finite algebra of global dimension $2$ with $r^2=0$.

 (2) $\tau S(3)\simeq S(4)$ and $\tau I(3)\simeq S(3)$. So $S(3)$ and $I(3)$ are $\tau$-rigid by Theorem 3.6.

 (3) By Theorem 3.12 or Theorem 3.15, Corollary 3.14 and (2), every indecomposable
 $\Lambda$-module is $\tau$-rigid.

 \vspace{0.2cm}

To show Theorem 3.12 and Theorem 3.15, in the following we will
construct an algebra $\Lambda$ and an indecomposable
$\Lambda$-module $M$ such that $\Lambda$ is of finite type and
finite global dimension with $r^2=0$ and there is no non-zero direct
summand of $P_0(M)$ and $P_1(M)$ in common.

\noindent{\bf Example 5.3} Let $\Lambda$ is given by the quiver:
 $$
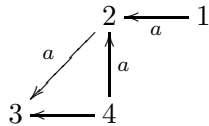
{\xymatrix{
 \ &2\ar[ld]_a& 1\ar[l]^a\\
 3&4\ar[u]_a\ar[l]&\ }}$$ with relations $a^2=0$. Then

(1) $\Lambda$ is a representation finite algebra of global dimension
$2$ with $r^2=0$.

(2) The injective module $I(3)$ has a minimal projective
resolution:\ \ \ \ \ \ \ \ \ \ \ \ \\
 \indent \ \ \ \  $P(2)\bigoplus P(3)\rightarrow P(2)\bigoplus P(4)\rightarrow
I(3)\rightarrow 0.$ By Theorem 3.12 (3), $I(3)$ is not $\tau$-rigid.

(3) $\tau^2 I(3)\simeq S(2)$, then by Theorem 3.6 (2) $\tau I(3)$ is
$\tau$-rigid. So $\tau M$ is $\tau$-rigid can not imply that $M$ is
$\tau$-rigid in general.

\vspace{0.3cm}

{\bf Acknowledgement } Part of the paper was written when the author
was visiting Bielefeld University in October and November of 2012
with the support of CRC 701. The author would like to thank Prof.
Henning Krause for invitation and hospitality. He also wants to
thank Prof. Claus Michael Ringel, Dr. Zhe Han and other people in
Bielefeld for useful discussion and kind help. The author is also
indebted to Prof. Osamu Iyama for kind help and Tiwei Zhao for
careful reading. The research of the author is carried out with the
support of NSFC (Nos.11101217, 11171142), NSF of Jiangsu Province
(No. BK20130983) and NSF for Colleges and Universities in Jiangsu
Province of China (No.11KJB110007).

\end{document}